\documentstyle[11pt,amscd,amssymb]{amsart}

\newcommand{\seq}[3]{{#1}_{#2}, \ldots, {#1}_{#3}}

\newcommand{\M}{M_{\Sigma}}
\newcommand{\Y}{\Sigma \times {\Bbb S}^1}
\newcommand{\X}{\Sigma \times {\Bbb P}^1}
\newcommand{\Spz}{\text{Sp}\, (2g,{\Bbb Z})}
\newcommand{\la}{\langle}
\newcommand{\ra}{\rangle}
\newcommand{\Dws}{D^{(w,\Sigma)}}
\newcommand{\ima}{{\bf i}}
\newcommand{\Cabg}{{\Bbb C}[\alpha,\beta,\gamma]}

\newcommand{\PGL}{\text{PGL}(2,{\Bbb C})}

\newcommand{\surj}{\twoheadrightarrow}

\newcommand{\ar}{\rightarrow}
\newcommand{\bd}{\partial}
\newcommand{\x}{\times}
\newcommand{\ox}{\otimes}
\newcommand{\iso}{\cong}
\newcommand{\isom}{\stackrel{\simeq}{\ar}}

\newcommand{\point}{\text{pt}}

\newcommand{\Hom}{\text{Hom}}
\newcommand{\Sym}{\text{Sym}}

\newcommand{\PD}{\text{P.D.}}

\newcommand{\ch}{\text{ch}\:}

\newcommand{\Todd}{\text{Todd}\,}

\newcommand{\cB}{{\cal B}}

\newcommand{\cE}{{\cal E}}

\newcommand{\cK}{{\cal K}}

\newcommand{\cO}{{\cal O}}

\newcommand{\cR}{{\cal R}}

\newcommand{\cU}{{\cal U}}
\newcommand{\cV}{{\cal V}}

\renewcommand{\AA}{{\Bbb A}}
\newcommand{\CC}{{\Bbb C}}

\newcommand{\EE}{{\Bbb E}}

\newcommand{\PP}{{\Bbb P}}

\newcommand{\SS}{{\Bbb S}}
\newcommand{\ZZ}{{\Bbb Z}}

\renewcommand{\a}{\alpha}
\renewcommand{\b}{\beta}

\newcommand{\g}{\gamma}
\newcommand{\e}{\varepsilon}

\renewcommand{\l}{\lambda}

\renewcommand{\o}{\omega}
\newcommand{\p}{\phi}
\newcommand{\q}{\psi}

\renewcommand{\S}{\Sigma}
\newcommand{\D}{\Delta}
\renewcommand{\L}{\Lambda}
\renewcommand{\P}{\Phi}
\newcommand{\Q}{\Psi}

\newcommand{\frg}{{\frak g}}

\theoremstyle{plain}
\begingroup
\theorembodyfont{\sl}
\newtheorem{thm}{Theorem}
\newtheorem{cor}[thm]{Corollary}
\newtheorem{lem}[thm]{Lemma}
\newtheorem{prop}[thm]{Proposition}

\endgroup

\theoremstyle{definition}

\theoremstyle{remark}

\newtheorem{ex}[thm]{Example}

\title{Gromov-Witten invariants of the moduli of bundles on a surface}

\author{Vicente Mu\~noz}
\address{Departamento de Matem\'aticas \\ Facultad de Ciencias \\
Universidad Aut\'onoma de Madrid \\ 28049 Madrid \\ Spain}
\email{vicente.munoz@@uam.es}
\thanks{Mathematics Subject Classification. Primary: 57R57. 
Secondary: 14D20, 58D27.}
\date{October, 1999}

\begin{document}

\maketitle

\section{Introduction}
\label{sec:1}

Let $\S=\S_g$ be a compact Riemann surface of genus $g \geq 2$
and let $\M$ stand for the moduli space of rank two stable vector 
bundles on $\S$ with odd (and fixed) determinant. In~\cite{quantum}
the author produced a presentation for the quantum cohomology
ring $QH^*(\M)$ in terms of its natural generators by giving
the relations satisfied by them. Here we want to show that this
information yields all the multiple-point Gromov-Witten invariants
on the generators, or which is equivalent, all $3$-point
Gromov-Witten invariants on the elements which are
quantum products of the generators.

On the other hand consider the (instanton) Floer cohomology
$HF^*(Y)$ of the three-manifold $Y=\Y$ endowed
with the $SO(3)$-bundle with
$w_2=\PD[\SS^1] \in H^2(Y;\ZZ/2\ZZ)$, which is determined
in~\cite{floer}. We show that the ring structure of
$HF^*(Y)$ yields the Donaldson invariants $D^{c_1}_{S}$ 
of the algebraic surface $S=\X$, for K\"ahler metrics whose 
period point is close enough to $\S$ in the 
K\"ahler cone, $U(2)$-bundles whose first Chern class $c_1
\in H^2(S,\ZZ)$ satisfies $c_1\cdot\S\equiv1\pmod 2$, and
on any collection of homology classes coming from $\S\subset S$.
Moreover the isomorphism $QH^*(\M)\iso HF^*(\Y)$ gives an
equality between these Donaldson invariants of $S$ and the
multiple-point Gromov-Witten invariants of $\M$ on the generators. 

What the quantum cohomology $QH^*(\M)$ does not give is the
$3$-point Gromov-Witten invariants on arbitrary homology classes
(i.e.\ on the elements which are cup products of the generators). 
This is equivalent to knowing the quantum product of two arbitrary
homology classes, and hence to knowing all the multiple-point 
Gromov-Witten invariants on homology classes. Again another 
equivalent formulation to this problem is to obtain the 
(just vector space) isomorphism
$$
  QH^*(\M) \stackrel{\iso}{\ar} H^*(\M).
$$
For the case $g=2$ this problem is settled in~\cite{Donaldson}.
Here we deal with the case $g=3$. For this it proves
necessary to write down all $3$-point Gromov-Witten invariants of 
degree $1$, which are computed thanks to the results 
in~\cite[section 3]{quantum} and moreover we need
two particular Gromov-Witten 
invariants of degree $2$, for which we use the results 
of~\cite{Sam}.

This may be viewed as an indication that (a presentation of) 
the quantum 
cohomology ring of a symplectic manifold has less information
than the full collection of Gromov-Witten invariants.

{\em Acknowledgements:\/} The author would like to thank 
Sharad Agnihotri and Marcos Mari\~no for asking him
to clarify these questions.

\section{Multiple-point Gromov-Witten invariants of $\M$}
\label{sec:2}

The general reference for this section is~\cite{quantum}.

\subsection*{Generators of the cohomology of $\M$}
Consider a symplectic basis $\{\seq{\g}{1}{2g}\}$ of 
$H_1(\S;\ZZ)$. Also $x\in H_0(\S;\ZZ)$ will stand for 
the class of the point. The $\mu$-map provide natural 
generators $\a=2\mu(\S) \in H^2(\M)$, $\q_i=\mu(\g_i) 
\in H^3(\M)$, $1\leq i \leq 2g$, and $\b=-4\mu(x) \in 
H^4(\M)$ of the cohomology ring $H^*(\M)$. We also put 
$$
  \AA(\S)=\Sym^*(H_0(\S)\oplus H_2(\S))\ox\L^* H_1(\S) 
  \stackrel{\mu}{\iso} \CC[\a,\b]\ox \L^*(\seq{\q}{1}{2g}),
$$
with grading $\deg\a=2$, $\deg \q_i=3$ and $\deg\b=4$. So a 
basic element of $\AA(\S)$ is of the form $z=\a^a\b^b
\q_{i_1}\cdots \q_{i_r}$, $a,b\geq 0$, $1\leq i_1<\cdots <i_r\leq
2g$. The cohomology class corresponding to $z$ under the
obvious epimorphism $\pi_H:\AA(\S)\surj H^*(\M)$ will be denoted
by 
$$
  \a_a\b_b\q_{i_1}\cup\cdots \cup\q_{i_r}= 
  \a\cup\stackrel{(a)}{\ldots}\cup \a \cup
  \b\cup\stackrel{(b)}{\ldots}\cup \b \cup \q_{i_1}
  \cup \cdots \cup \q_{i_r},
$$
and leave the notation with exponents for the quantum
product (note that we cannot suppress the subindices
$a$ and/or $b$ when they are $1$).

The action of the mapping class group of $\S$ yields 
an action of the symplectic group $\Spz$ on $\{\q_i\}$ 
and the invariant part $H^*(\M)_I$ of $H^*(\M)$ 
is generated by $\a$, 
$\b$ and $\g=-2 \sum_{i=1}^g \q_i\cup\q_{g+i}$.
Also we denote $\a_a\b_b\g_c=
  \a\cup\stackrel{(a)}{\ldots}\cup \a \cup
  \b\cup\stackrel{(b)}{\ldots}\cup \b \cup
  \g\cup\stackrel{(c)}{\ldots}\cup \g$ for a typical
cohomology class in $H^*(\M)_I$.

\subsection*{Gromov-Witten invariants of $\M$}
Let $A$ denote the positive generator of $\pi_2(\M)
\iso H_2(\M;\ZZ)$, i.e.\ $\a[A]=1$. Fix $d\geq 0$ and 
$r\geq 3$ and let $z_i\in H^{p_i}(\M)$, $1\leq i\leq r$,
be homogeneous 
homology classes. The $r$-point Gromov-Witten invariant
$\Q_{dA}^{\M}(z_1, \ldots ,z_r)$ is defined to be zero 
if $p_1+\cdots+p_r\neq 6g-6+4d$, and in the following
manner if $p_1+\cdots+p_r=6g-6+4d$. Consider $r$ 
different points $\seq{P}{1}{r} \in {\PP}^1$ and 
represent the Poincar\'e duals of $z_i$ by generic 
cycles $V_{z_i}\subset \M$. Then 
$$
  \Q_{dA}^{\M}(\seq{z}{1}{r}) = \#\{f:\PP^1 \ar \M | 
  \text{$f$ is holomorphic}, f_*[\PP^1]=dA,
  f(P_i)\in V_{z_i}, 1\leq i\leq r \} 
$$
where $\#$ denotes count with signs. This is well-defined 
in appropriate circumstances~\cite{quantum}.
We also put $\Q^{\M}=\sum_{d\geq0}\Q_{dA}^{\M}$ and
extend the definition to non-homogeneous $z_i$ multi-linearly.

\subsection*{Quantum cohomology of $\M$}
The quantum cohomology $QH^*(\M)$ is $H^*(\M)$ as a vector space
(so they are identified by the identity map), but
the ring structure of $QH^*(\M)$, the quantum multiplication, 
is a deformation of the usual cup product for $H^*(\M)$, and it
is graded only modulo $4$. It is defined as follows.
Let $r \geq 2$ be any integer. For cohomology classes $z_i 
\in H^*(\M)$, $1 \leq i \leq r$, the quantum product of these
classes is $z_1 \cdots z_r \in H^*(\M)$ defined by 
$$
  \la z_1 \cdots z_r, z_{r+1} \ra = 
  \Q^{\M}(z_1, \ldots ,z_r, z_{r+1}), \text{ for any }
   z_{r+1} \in H^*(\M).
$$

Obviously, the $3$-point Gromov-Witten invariant completely 
determines the quantum product. Therefore the $3$-point 
Gromov-Witten invariant determines the multiple-point 
Gromov-Witten invariant by the simple formula
$$
  \Q^{\M}(z_1,\ldots,z_r)= \Q^{\M}(z_1,z_2,z_3\cdots z_r),
$$
for any $z_1,\ldots,z_r\in H^*(\M)$ homology classes, $r\geq 3$.

The ring $QH^*(\M)$ is generated by $\a$, $\b$ and $\q_i$,
$1\leq i\leq 2g$. So there is an epimorphism 
$\pi_{QH}: \AA(\S)\surj QH^*(\M)$.
The typical quantum product is denoted as
$$
  \a^a\b^b \q_{i_1}\cdots \q_{i_r} = 
  \a\stackrel{(a)}{\cdots} \a
  \b\stackrel{(b)}{\cdots} \b \q_{i_1}\cdots \q_{i_r}.
$$

\subsection*{Multiple-point Gromov-Witten invariant on generators}
We define the {\em multiple-point Gromov-Witten invariant of $\M$
on generators\/} in the following way.
For any $z=\a^a\b^b\q_{i_1}\cdots\q_{i_r}\in \AA(\S)$, put
$$
  \tilde{\Q}^{\M}(z)= \Q^{\M}(\a, \stackrel{(a)}{\ldots}, \a, \b ,
  \stackrel{(b)}{\ldots} ,\b,\q_{i_1}, \ldots, \q_{i_r}),
$$
then extend by linearity. Another way to say this is
for $z\in \AA(\S)$, 
$\tilde{\Q}^{\M}(z)=\la \pi_{QH}\,z\ra_{\M}$, where
$\la \cdot\ra_{\M}:H^*(\M) \ar \CC$ is the
pairing with the volume class. Note that 
$\tilde\Q^{\M}$ is invariant under the action of $\Spz$, so it is
determined by its effect on invariant elements $z\in \AA(\S)_I$.
For $g\geq 3$, $\g=-2\sum_{i=1}^g \q_i\q_{g+i} \in QH^*(\M)$, 
by~\cite[lemma 14]{quantum}, so for $z=\a^a\b^b\g^c\in \AA(\S)$ 
we have 
$$
  \tilde\Q^{\M}(z)= \Q^{\M}(\a,\stackrel{(a)}{\ldots},\a,\b, 
  \stackrel{(b)}{\ldots},\b,\g, \stackrel{(a)}{\ldots},\g).
$$

  $\tilde\Q^{\M}$ determines the quantum cohomology $QH^*(\M)$.
  Actually $QH^*(\M)=\AA(\S)/I$, for an ideal $I\subset \AA(\S)$
  of relations. Then $I=\{R\in \AA(\S)| \tilde\Q^{\M}(Rz)=0, 
  \forall z\in \AA(\S) \}$. The following result is a converse
  of this.

\begin{prop}
\label{prop:1}
  The quantum cohomology of $\M$ and the value of $\la
  \g_{g-1}\ra_{\M}= 2^{g-1}g!$ determines 
  the multiple-point Gromov-Witten invariant on generators
  $\tilde\Q^{\M}$.
\end{prop}

\begin{pf}
  Let $z\in \AA(\S)$. As $\tilde\Q^{\M}(z)$ is invariant 
  under the action of $\Spz$, we may project $z$ to the
  invariant part $\AA(\S)_I\subset \AA(\S)$, i.e.\ we may
  suppose that $z$ is invariant. A basis for $QH^*(\M)_I$ is 
  $\a^a\b^b\g^c$, $a+b+c<g$, by~\cite[section 5]{quantum}. 
  This means that
  $z=R+\sum\limits_{a+b+c<g}c_{abc}\a^a\b^b\g^c$, 
  for some numbers
  $c_{abc}\in \CC$ and a relation $R\in \AA(\S)_I$ for $QH^*(\M)_I$.
  By definition of relation we have that $\tilde\Q^{\M}(R)=0$. Also 
  $\tilde\Q^{\M}(\a^a\b^b\g^c)=0$, for $a+b+c<g$ and
  $(a,b,c) \neq (0,0,g-1)$, 
  since in that case $\deg(\a^a\b^b\g^c) =2a+4b+6c< 6g-6$. Hence
  $\tilde\Q^{\M} (z)= c_{0,0,g-1} 
  \tilde\Q^{\M} (\g^{g-1})$.
  Finally the value of $\tilde\Q^{\M} (\g^{g-1})=
  \la \g_{g-1},[\M]\ra=\la \g_{g-1}\ra_{\M}= 2^{g-1} g!$,
  by~\cite{Thaddeus}.
\end{pf}

Let $z_i=\a^{a_i}\b^{b_i}\seq{\q}{i,j_1}{i,j_{r_i}}$, $i\in I$, 
be a collection of elements of $\AA(\S)$ such that $\cB_H=
\{\pi_H z_i\}_{i\in I}$ is a basis for $H^*(\M)$. 
Then $\cB_{QH}=\{\pi_{QH} z_i\}_{i\in I}$ will be
a basis for $QH^*(\M)$. The following statements are equivalent
\begin{enumerate}
\item We have a presentation of $QH^*(\M)$, i.e.\ we know the ideal 
  $I\subset \AA(\S)$ such that $QH^*(\M)=\AA(\S)/I$.
\item We know how to compute $\tilde\Q^{\M}(z)$, for any $z\in \AA(\S)$.
\item We know the  $3$-point Gromov-Witten invariant on quantum products
of the generators, i.e.\ $\Q^{\M}(\pi_{QH}z_i,\pi_{QH}z_j,\pi_{QH}z_k)$,
for all $i,j,k\in I$.
\item We know the coefficients $c_{ijk}$ such that $\pi_{QH}z_i
  \,\pi_{QH}z_j=\sum_k c_{ijk} \pi_{QH}z_k$, for $i,j\in I$.
\end{enumerate}
The equivalence of 1 and 2 is proposition~\ref{prop:1}. The equivalence of
2 and 3 follows from $\Q^{\M}(\pi_{QH}z_i,\pi_{QH}z_j,\pi_{QH}z_k)=
\tilde\Q^{\M}(z_iz_jz_k)$. The equivalence of 3 and 4 follows from the 
definition (using the intersection pairing). The statement 1 is true 
because of~\cite{quantum}, so the other three follow, meaning at least
that we can perform all computations for any fixed genus $g$ by hand or with
a computer.

A full knowledge of the Gromov-Witten invariants requires
to know the $3$-point Gromov-Witten invariants on homology
classes $\pi_H z_i=\a_{a_i}\b_{b_i}
\q_{i,j_1}\cup\cdots\cup \q_{i,j_{r_i}}$, $i\in I$. The following 
statements are equivalent
\begin{enumerate}
\item We know the  $3$-point Gromov-Witten invariant on arbitrary
homology classes, i.e.\ $\Q^{\M}(\pi_Hz_i,\pi_Hz_j,\pi_Hz_k)$, for
any $i,j,k\in I$.
\item We know the coefficients $d_{ijk}$ such that $\pi_{H}z_i
  \,\pi_{H}z_j=\sum_k d_{ijk} \pi_{H}z_k$, for any $i,j\in I$.
\item We know the isomorphism $H^*(\M) \isom QH^*(\M)$ in terms of the
basis $\cB_H$ and $\cB_{QH}$, respectively.
\end{enumerate}
The equivalence of 1 and 2 is obvious using the intersection pairing.
If we have 3 then 2 follows by translating the point 4 above (which
we already have) through the isomorphism. Conversely, if we have 2
we may quantum multiply the generators 
$\a$, $\b$ and $\q_i$, $1\leq i\leq 2g$, repeteadly until we get
any $\pi_{QH}z_i$ in terms of the basis $\cB_H$.
All this information is stronger than a presentation of the
quantum cohomology ring, as we shall see in the example of
section~\ref{sec:4}.

\section{Donaldson invariants of $S=\X$}
\label{sec:3}

The general reference for this section is~\cite{floer}.

\subsection*{Donaldson invariants of $S$}
Consider the algebraic surface $S=\X$. In particular $S$ is 
a smooth $4$-manifold with $b_1=2g$ and $b^+=1$. The Donaldson 
invariants of $S$ depend on a metric as explained in~\cite{wall}.
The K\"ahler cone of $S$ is $\cK=\{a\S+b\PP^1 | a,b>0\} \subset H^2(S)$. 
We set $\AA(S)=\Sym^*(H_0(S) \oplus H_2(S))
\ox \L^* H_1(S)$, which is graded giving degree $4-i$ to the
classes in $H_i(S)$. The class of the point will be denoted
by $x\in H_0(S)$. The Donaldson invariants are linear functionals
$$
  D^{c_1}_{S,H}:\AA(S) \ar \CC,
$$
depending on the first Chern class $c_1 \in H^2(S;\ZZ)$ 
of the $U(2)$-bundles involved, and on a (generic) polarisation 
$H\in\cK$. 

We shall consider only polarisations close to $\S$, i.e.\
$H=\S+\e \PP^1$, with $\e>0$ small (how small depending on
the degree of the element on which we are computing the
Donaldson invariants). So $D^{c_1}_S$ shall stand for 
$D^{c_1}_{S,H}$ with such polarisation $H$. Also we only 
consider the cases with $c_1\cdot \S$ odd. As
$D^{c_1+2\a}_S=(-1)^{\a^2}D^{c_1}_S=D^{c_1}_S$, for $\a
\in H^2(S;\ZZ)$, we may suppose that $c_1=\PP^1$ or $c_1=
\PP^1+\S$. If $c_1=\PP^1$ then $D^{c_1}_S$ is non-zero 
only on elements of degree $6g-6+4d$, $d\geq 0$ even. If
$c_1=\PP^1+\S$ then $D^{c_1}_S$ is non-zero only on 
elements of degree $6g-6+4d$, $d\geq 0$ odd. 
We collect all these invariants together by putting
$w=\PP^1$ and $\Dws_S=D^w_S+D^{w+\S}_S$.

A basic element $z\in \AA(S)$ is of the form
$z=\S^a x^b\g_{i_1}\cdots\g_{i_r} (\PP^1)^{e}$. We restrict
to classes coming from $\S \subset S$, i.e.\ $z\in \AA(\S)
\subset \AA(S)$. Therefore $e=0$ and $z=\S^a x^b\g_{i_1}
\cdots\g_{i_r}$. As in section~\ref{sec:2} we identify $\AA(\S)\iso 
\CC[\a,\b]\ox \L^*(\seq{\q}{1}{2g})$ by $\a=2\S$, 
$\b=-4x$ and $\q_i=\g_i$, $1\leq i\leq 2g$. Also we let
$\g=-2 \sum \q_i\q_{g+i}$.

\subsection*{Floer cohomology of $\Y$}
Consider the three-manifold $Y=\Y$ and let $HF^*(Y)$ stand
for its (instanton) Floer cohomology endowed with the 
$SO(3)$-bundle with $w_2=\PD[\SS^1]$. By the analysis in~\cite{floer},
$HF^*(Y)$ is generated by $\a$, $\b$, $\q_i$, $1\leq i\leq 2g$, coming
from the relative Donaldson invariants for $\S\x D^2$
of \, $2\S$, $-4x$, $\g_i$, $1\leq i\leq 2g$, respectively.

There is an evaluation map~\cite[corollary 19]{floer}
$\la \cdot \ra_g: HF^*(Y)\ar \CC$ with the property that
for any $z=\a^a\b^b\q_{i_1}\cdots\q_{i_r} \in \AA(\S)$, we have
$$
  \Dws_S (z)= \Dws_S ((2\S)^a(-4 x)^b \g_{i_1}
  \cdots\g_{i_r}) = \la \a^a\b^b\q_{i_1} \cdots\q_{i_r} \ra_g,
$$

The Donaldson invariants $\Dws_S$ determine the Floer homology since
$HF^*(Y)=\AA(\S)/J$ and
$J=\{R\in \AA(\S) | \Dws_S(Rz)=0, \forall z\in\AA(\S)\}$.
The following result is a converse of this.

\begin{prop}
\label{prop:2}
  The Floer cohomology of $Y=\Y$ and the value $D_S^w (\g^{g-1})
  =- 2^{g-1}g!$ determine the Donaldson invariants $\Dws_S (z)$, 
  for any $z\in \AA(\S)$.
\end{prop}

\begin{pf}
  The argument is similar to that of proposition~\ref{prop:1}. 
  We may take $z\in \AA(\S)$ invariant under $\Spz$. A vector
  basis of $HF^*(Y)_I$ is $\a^a\b^b\g^c$, $a+b+c<g$,
  by~\cite{floer}. Then $z=R+\sum\limits_{a+b+c<g}c_{abc}
  \a^a\b^b\g^c$, for some $c_{abc}\in \CC$ and $R\in \AA(\S)_I$ 
  which is a relation for $HF^*(Y)_I$. Then $\Dws_S(R)=0$ and 
  $\Dws_S(\a^a\b^b\g^c)=0$, for $a+b+c<g$ and  $(a,b,c) \neq 
  (0,0,g-1)$, since in that case $\deg(\a^a\b^b\g^c) =2a+4b+6c< 
  6g-6$ and there is no moduli space of ASD connections of
  dimension smaller than $6g-6$. Hence
  $\Dws_S (z)= c_{0,0,g-1} D^w_S (\g^{g-1})= -c_{0,0,g-1}
  \la \g_{g-1},[\M]\ra$.
\end{pf}

\begin{thm}
\label{thm:3}
  Suppose that $g\geq 3$. Then for any 
  $z=\a^a\b^b\q_{i_1}\cdots\q_{i_r}\in \AA(\S)$ 
  with $\deg z=6g-6+4d$, $d\geq 0$, we have 
$$
  \Q^{\M}_{dA}(\a, \stackrel{(a)}{\ldots}, \a, \b ,
  \stackrel{(b)}{\ldots} ,\b,\q_{i_1}, \ldots, \q_{i_r})=(-1)^{gd+1} 
  \Dws_S ((2\S)^a(-4x)^b \g_{i_1} \cdots\g_{i_r}).
$$
\end{thm}

\begin{pf}
  By~\cite[corollary 21]{quantum}, for any $g \geq 3$ there is 
  an isomorphism $QH^*(\M) \isom HF^*(Y)$, $(\a,\b,\q_i) \mapsto 
  (\e^{2g}\a,\e^{4g}\b,\e^{3g}\q_i)$, where $\e$ is a primitive 
  eighth root of unity. As $D_S^w (\g^{g-1})=-\tilde\Q^{\M} 
  (\g^{g-1})$ the result follows easily.
\end{pf}

\begin{cor}
  Suppose that $g\geq 4$. Then for any 
  $z=\a^a\b^b\q_{i_1}\cdots\q_{i_r}\in \AA(\S)$ 
  with $\deg z=6(g-1)-6+4d$, $d\geq 0$, we have   
  $$
  \Q^{M_{\S_g}}_{dA}(\g, \a, \stackrel{(a)}{\ldots}, \a, \b ,
  \stackrel{(b)}{\ldots} ,\b,\q_{i_1}, \ldots, \q_{i_r}) =
 (-1)^d 2g \Q^{M_{\S_{g-1}}}_{dA}(\a, \stackrel{(a)}{\ldots}, 
  \a, \b ,
  \stackrel{(b)}{\ldots} ,\b,\q_{i_1}, \ldots, \q_{i_r}).
  $$
\end{cor}

\begin{pf}
  By~\cite[corollary 19]{floer} it is 
  $\Dws_{\S_g\x\PP^1}(\g z)= 2g
  \Dws_{\S_{g-1}\x\PP^1}(z)$, for any $z\in\AA(\S)$. 
  Then the result follows from theorem~\ref{thm:3}.
\end{pf}

\subsection*{Generating function}
Proposition~\ref{prop:2} can be used for a given genus $g$ to 
effectively compute the Donaldson invariants $\Dws_S$ on homology
classes coming from $\S\subset S$. For this, we collect
the Donaldson invariants into a generating function as follows 
(we only need the invariant part)
\begin{equation}
\label{eqn:gener}
  \Dws_S(e^{s\a+\l \b+r\g})=\la e^{s\a+\l \b+r\g} \ra_g.
\end{equation}
Every relation $R(\a,\b,\g)$ for $HF^*(Y)_I$ gives a differential
equation $R({\bd\over\bd s},{\bd\over \bd\l},{\bd\over\bd r})$ 
satisfied by~\eqref{eqn:gener}. Actually~\eqref{eqn:gener} is the
only solution $F(s,\l,r)$ to all these differential equations satisfying the
initial conditions
$$
\left\{ \begin{array}{ll}
  {\bd^{a+b+c}F\over\bd s^a\bd \l^b\bd r^c}\big|_{s=\l=r=0}=0, & a+b+c<g, 
  \  (a,b,c)\neq (0,0,g-1) \\
  {\bd^{g-1}F\over \bd r^{g-1}}\big|_{s=\l=r=0} =-2^{g-1} g!
\end{array} \right.
$$

The decomposition $HF^*(Y)_I=\bigoplus\limits_{r=-(g-1)}^{g-1} 
R_{g,r}$ in~\cite[section 7]{floer} may be useful
for finding $F(s,\l,r)$. The map $\la \cdot \ra_g:
HF^*(Y)_I \ar \CC$ yields maps $\la \cdot\ra_{g,r}:R_{g,r}\ar \CC$ 
such that $\la \cdot\ra_g=\sum\limits_{r=-(g-1)}^{g-1}\la \cdot\ra_{g,r}$. 
Now let $\Phi_r=\la e^{s\a+\l \b+r\g} \ra_{g,r}$. Then
the relations for $R_{g,r}$ give the differential equations
satisfied by every $\Phi_r$ and $\Dws_S(e^{s\a+\l \b+r\g})=\sum_r
\Phi_r$.

\begin{ex}
\label{ex:g=2}
  Let the genus $g=2$. By~\cite[example 22]{floer}
  $$
  HF^*(Y)_I=R_{2,-1}\oplus R_{2,0}\oplus R_{2,1}=
  {\Cabg \over (\a-4, \b+8, \g)} \oplus {\Cabg \over 
  (\a^2, \b- 8,\g+16\a)} \oplus {\Cabg \over (\a+4, \b+8, \g)}.
  $$
  Associated to $R_{2,-1}$ we have the differential equations 
  ${\bd \over\bd s}-4$, ${\bd \over\bd \l}+8$ and 
  ${\bd \over\bd r}$, whose solution is a multiple of
  $e^{4s-8\l}$.
  For $R_{2,0}$ we have the differential equations 
  ${\bd^2 \over\bd s^2}$, ${\bd \over\bd \l}-8$ and 
  ${\bd \over\bd r}+16{\bd \over\bd s}$, 
  whose solutions are linear combinations of $e^{8\l}$
  and $(16r-s)e^{8\l}$. Finally for $R_{2,1}$ we get
  $e^{-4s-8\l}$. Then $\Dws_S(e^{s\a+\l \b+r\g})$
  is a linear combination of $e^{4s-8\l}$, $e^{8\l}$,
  $(16r-s)e^{8\l}$ and $e^{-4s-8\l}$. 
  Putting the initial conditions one gets that
$$
\Dws_S(e^{s\a+\l \b+r\g})=
-\frac{1}{16} \sinh 4s\, e^{-8\l} -\frac14 (16r-s)  e^{8\l}.
$$
\end{ex}

\begin{ex}
\label{ex:g=3}
  In the case of genus $g=3$ we use the decomposition 
  in~\cite[example 23]{floer} to obtain that 
$$
\Dws_S(e^{s\a+\l \b+r\g})=
\frac{1}{2048} \cos 8s \,e^{8\l} +\frac{3}{128} \cosh 4s\, 
e^{-8\l}  
$$
$$
-\frac{1}{32}(s\sinh 4s-4\l \cosh 4s +12r\sinh 4s)
e^{-8\l} - \left(\frac{49}{2048} -\frac14 \l+\frac{3}{64}
(16r-s)^2 +\l^2\right) e^{8\l}
$$
  This agrees with~\cite[formula (5.16)]{Marino}.
\end{ex}

Examples~\ref{ex:g=2} and~\ref{ex:g=3} agree with the fact that
${\bd\over\bd r}\Dws_{\S_g\x\PP^1}(e^{s\a+\l \b+r\g})=\la \g e^{s\a+\l \b+r\g} 
\ra_g=2g \la e^{s\a+\l \b+r\g} \ra_{g-1}=2g
\Dws_{\S_{g-1}\x\PP^1}(e^{s\a+\l \b+r\g})$. Note
also that by theorem~\ref{thm:3}, 
$$
 \tilde{\Q}^{\M}(e^{s\a+\l \b+r\g})=\left\{ 
  \begin{array}{ll} (-1)^{g+1\over 2}\Dws_S(e^{\ima s\a-\l \b-
  \ima r\g}) \qquad & \text{$g$ odd} \\  
 -\Dws_S(e^{s\a+\l \b+r\g}) & \text{$g$ even} \end{array}\right.
$$

\section{Gromov-Witten invariants in the case $g=3$}
\label{sec:4}

In this section we shall compute all the Gromov-Witten invariants 
for the moduli space
$\M$ of rank $2$ odd degree (and fixed determinant) stable bundles over
the Riemann surface of genus $g=3$. The dimension of $\M$ is
$\dim\M=6g-6=12$ and its cohomology ring is written in~\cite[proposition 
1]{quantum}
to be
\begin{equation}
H^*(\M)={\Cabg\over(q^1_3,q^2_3,q^3_3)} 
  \oplus \left( H^3 \ox {\Cabg\over(q^1_2,q^2_2,q^3_2)}\right)
  \oplus
 \left( \L^2_0 H^3 \ox {\Cabg\over(q^1_1,q^2_1,q^3_1)}\right), 
\label{eqn:1}
\end{equation}
where $H^3=\la \seq{\q}{1}{6}\ra \iso H_1(\S)$, 
$\L^2_0 H^3=\la \q_i\cup\q_j | 1\leq i<j\leq 6, j\neq i+3\ra \oplus
\la\q_1\cup\q_4-\q_2\cup\q_5,\q_1\cup\q_4-\q_3\cup\q_6\ra$ and
$$
\begin{array}{lll}
  q^1_1=\a,  & q^1_2=\a_2 +\b, \qquad & q^1_3=\a_3+5\a_1\b_1 +4\g,  \\
  q^2_1=\b,  & q^2_2=\a_1\b_1 +\g, & q^2_3=\a_2\b_1+\b_2 +\frac43\g_1\a_1,\\
  q^3_1=\g,\qquad & q^3_2=\g_1\a_1, & q^3_3=\g_1\a_2+\g_1\b_1.
\end{array}
$$
The quantum cohomology is~\cite[theorem 20]{quantum}
\begin{equation}
QH^*(\M)={\Cabg\over(Q^1_3,Q^2_3,Q^3_3)} 
  \oplus \left( H^3 \ox {\Cabg\over(Q^1_2,Q^2_2,Q^3_2)}\right)
  \oplus
 \left( \L^2_0 H^3 \ox {\Cabg\over(Q^1_1,Q^2_1,Q^3_1)}\right), 
\label{eqn:2}
\end{equation}
where $\L^2_0 H^3=\la \q_i\q_j | 1\leq i<j\leq 6, j\neq i+3\ra \oplus
\la\q_1\q_4-\q_2\q_5,\q_1\q_4-\q_3\q_6\ra$ and
$$
\begin{array}{lll}
  Q^1_1=\a,  & Q^1_2=\a^2 +\b+8,\qquad & Q^1_3=\a^3+5\a\b +4\g-24\a,\\
  Q^2_1=\b+8, \qquad & Q^2_2=\a\b+\g-8\a, & Q^2_3=\a^2\b+\b^2 
  +\frac43\g\a+8\a^2+16\b+64,\\
  Q^3_1=\g, & Q^3_2=\g\a, & Q^3_3=\g\a^2+\g\b+8\g.
\end{array}
$$
We shall determine the isomorphism $H^*(\M)\iso QH^*(\M)$, which is
equivalent to the Gromov-Witten invariants by the discussion at
the end of section~\ref{sec:2}. We remark that the three pieces of
the decompositions~\eqref{eqn:1} and~\eqref{eqn:2} correspond by 
the proof of~\cite[proposition 16]{quantum}.

\subsection*{The invariant part}
A basis for $H^*(\M)_I$ is $\cB_H=\{1,\a,\b,\a_2,\g,\a_1\b_1,\b_2,
\a_1\g_1,\b_1\g_1,\g_2\}$ and a basis for $QH^*(\M)_I$ is 
$\cB_{QH}=\{1,\a,\b,\a^2,\g,\a\b,\b^2,\a\g,\b\g,\g^2\}$.
The isomorphism $H^*(\M)_I \iso QH^*(\M)_I$ is given by
$$
\left\{\begin{array}{lclclclcl}
\a^2 &=& \a_2 &+& A_1 \\
\a\b &=&\a_1\b_1 &+&A_2 \a  \\
\b^2 &=&\b_2 &+&A_3 \b+A_4\a_2 &+&B_1 \\
\a\g &=&\a_1\g_1 &+&A_5\b+A_6 \a_2 &+&B_2 \\
\b\g &=&\b_1\g_1&+& A_7\g+A_8\a_1\b_1 &+&B_3 \a \\
\g^2 &=&\g_2 &+& A_9\g_1\a_1+A_{10}\b_2 &+&B_4 \a_2+B_5\b &+&C
\end{array}\right.
$$
for constants $\seq{A}{1}{10}$, $\seq{B}{1}{5}$ and $C$, determined 
by the Gromov-Witten invariants of degree $1$, $2$ and $3$, 
respectively. Now we shall compute these constants. To start with
we write the top-products from~\cite{Thaddeus}
(alternatively we may use $\la\g_2\ra_{\M}=24$ and the relations
for $H^*(\M)_I$)
$$
\begin{array}{llll}
  \la\a_6\ra_{\M}=7\cdot 32, \qquad & \la\a_4\b_1\ra_{\M}= -64,\qquad
   & \la\a_2\b_2\ra_{\M}= 32,\qquad &\la\b_3\ra_{\M}=0, \\ 
  \la\a_3\g_1\ra_{\M}=24, & \la\a_1\b_1\g_1\ra_{\M}=-24 & 
  \text{and}\quad \la\g_2\ra_{\M}=24.
\end{array}
$$

\subsection*{The Gromov-Witten invariants of degree one}
Let us recall the space $N$ given in~\cite[section 3]{quantum}, 
which parametrises non-split extensions on $\S$ of the form
$$
 0\ar L\ar E\ar \L\ox L^{-1}\ar 0,
$$
where $L$ is a line bundle of degree zero. 
Let $J$ be the Jacobian variety parametrising line
bundles of degree zero on $\S$. The natural isomorphism
$H_1(\S)\iso H^1(J)$ associates to $\{\g_i\}$ a symplectic
basis $\{\p_i\}$ of $H^1(J)$. Put $\o=\sum_{i=1}^3 \p_i\wedge
\p_{3+i}$. Then $N$ sits as a fibration $\PP^2\ar N=
\PP(\cE^{\vee})\ar J$, where $\cE$ is a rank-$3$ bundle 
over $N$ with Chern classes $c_i=c_i(\cE)=(4^i/i!)\o^i$.
The cohomology ring of $N$ is generated by the classes 
$\p_i \in H^1(N)$, $1\leq i\leq 6$ and the class of
the hyperplane $h\in H^2(N)$. Actually $H^*(N)$ is 
freely generated (in a graded sense) subject to the only
relation
$$
   h^3+c_1h^2+c_2h+c_3=0.
$$
(We only use cup-products in $H^*(N)$, so we write $\p_{i_1}\cdots
\p_{i_r}h^j=\p_{i_1}\cup\cdots\cup\p_{i_r}\cup h\cup
\stackrel{(j)}{\ldots}\cup h$.)
There is an action of the mapping class group of $\S$ on $N$, which
gives an action of $\Spz$ on $\{\p_i\}$ and the invariant part
is generated by $\o$ and $h$. 

There is a natural map $i:N\ar \M$. We need to write down
how our basis for $H^*(\M)_I$ restrict to $H^*(N)$. 
By~\cite[equation (12)]{quantum} $\a=4\o+h$, $\b=h^2$, 
$\q_i=-h\p_i$, $1\leq i\leq 6$, and $\g=-2\o h^2$ in $N$
(we omit $i^*$ as there is no danger of confussion).
Therefore
$$
\left\{\begin{array}{l}
  \a=4\o+h \\
  \b =h^2\\
  \a_2=16\o^2+8\o h+h^2 \\
  \g= -2\o h^2\\
  \a_1\b_1=4\o h^2+h^3= - 8\o^2 h -\frac{32}{3}\o^3 \\
  \b_2 =h^4= 8\o^2 h^2+\frac{64}{3} \o^3 h\\
  \a_1\g_1 = -8\o^2 h^2 -2\o h^3= 16\o^3 h\\
  \b_1\g_1 =-2\o h^3= -16 \o^3 h^2 \\
  \g_2= 0
\end{array}\right.
$$
Let $l\in H_2(N;\ZZ)$ be the class of the line in a fibre 
$\PP^2\subset N$. Then $i_*l=A$. The Gromov-Witten invariants of 
degree $1$ of $N$ and $\M$ satisfy $\P_l^N=\P_A^{\M}$ 
by~\cite[lemma 8]{quantum}. Moreover the Gromov-Witten invariants
of $N$ are computed as in~\cite[lemma 10]{quantum}. Let 
$z_i=\o^{a_i} h^{b_i} \in H^*(N)$, $1\leq i\leq 3$, $0\leq a_i
\leq 3$, $0\leq b_i\leq 2$. Put 
$a=a_1+a_2+a_3$ and $b=b_1+b_2+b_3$ and suppose $2a+2b=6g-2=16$. 
Then 
$$
\Q^N_l(z_1,z_2,z_3)=\left\{\begin{array}{ll}
 \la \o^a X^b,[J]\ra=
{(-8)^{b-5}\over (b-5)!}\o^3=6 {(-8)^{b-5}\over (b-5)!}, \qquad &
  b\geq 5 \\ 0, & b<5 \end{array} \right.
$$
where $X^{2g-1+i}=X^{5+i}=((-8)^i/i!)\o^i$ as 
in~\cite[lemma 10]{quantum}. With this we may find
the coefficients $\seq{A}{1}{10}$. For instance we have
$$
\left.
\begin{array}{l}
\Q_A^{\M}(\a,\a,\g_2)=\la \a^2,\g_2\ra =\la A_1,\g_2\ra=24 A_1 \\
\Q_l^N(\a,\a,\g_2)=\Q_l^N(\o+h,\o+h,0)=0
\end{array} \right\} \implies A_1=0,
$$
$$
\left.
\begin{array}{l}
\Q_A^{\M}(\a,\b,\b_1\g_1)=\la \a\b,\b_1\g_1\ra=\la A_2\a,
  \b_1\g_1\ra=-24 A_2  \\
\Q_l^N(\a,\b,\b_1\g_1)=\Q_l^N(4\o+h,h^2,-16\o^3h^2)=-16\cdot 6
\end{array} \right\} \implies A_2=4.
$$
Analogously we get $A_3=-12$, $A_4=-8$, 
$A_5=-3$, $A_6=-3$, $A_7=-20$, $A_8=-12$, 
$A_9=8$ and $A_{10}=-6$. 

\subsection*{More constants from pairings}
Once we know the coefficients $A_1,\ldots, A_{10}$, we may 
get some relations between the coefficients $B_1,\ldots,B_5$ and $C$
by working out the pairings of the
elements of the basis $\cB_{QH}$ written in terms of
those of $\cB_H$. The only pairings in which the 
coefficients $B_1,\ldots,B_5$ and $C$ are going to appear
are those between elements whose degrees add up to 
$6g-6+4d$, $d\geq 2$, i.e.\ $\la \b\g, \b\g\ra$,
$\la \g^2, \b^2\ra$, $\la \g^2, \g\a\ra$ 
and $\la \g^2, \g^2\ra$. The relations in $QH^*(\M)_I$ imply
that $\g^3=0$ and $\g^2(\b^2-64)=0$, so $\la \b\g, \b\g\ra=
\la \b^2\g^2\ra_{\M}=\la 64\g^2 \ra_{\M}=64\cdot 24$,
$\la \g^2, \b^2\ra=64\cdot 24$, $\la \g^2, \g\a\ra=0$ 
and $\la \g^2, \g^2\ra=0$. So we get
$$
  64\cdot 24= \la \b\g, \b\g\ra=
  \la \b_1\g_1 -20 \g -12 \a_1\b_1 +B_3 \a, 
  \b_1\g_1 -20 \g -12 \a_1\b_1 +B_3 \a \ra = \qquad
$$
$$ 
 \qquad = 2B_3 \la \a_1\b_1\g_1
  \ra_{\M} +\la -20 \g -12 \a_1\b_1,-20 \g -12 \a_1\b_1\ra=
  -24 \cdot 2B_3 +24\cdot 112,
$$
and hence $B_3=24$. From the other pairings we obtain the equations
$B_2+B_4-B_5=-24$, $B_1+\frac43 B_4=-32$ and 
$C=8B_5$, respectively. No more information can be extracted 
from the intersection pairing.

To finish we only need to find
$$
\left\{\begin{array}{l} 24B_1=\Q_{2A}^{\M} (\b,\b,\point) \\
 24B_2=\Q_{2A}^{\M} (\a,\g,\point) \end{array}\right.
$$
This correspond to the computation of two particular Gromov-Witten
invariants of degree two, task that will be carried out
in section~\ref{sec:degree-2}. The answer is given in 
equations~\eqref{eqn:ag}
and~\eqref{eqn:bb}. We get $B_1=0$, $B_2=-1$ and then
$B_5=-1$, $B_4=-24$ and 
$C=-8$.

\subsection*{Non-invariant part}
We recall that the decompositions~\eqref{eqn:1} and~\eqref{eqn:2} correspond 
since every piece is the isogeneous piece corresponding to an irreducible
representation of $\Spz$. First we shall deal with the piece corresponding to
$H^3$. A vector basis given by the usual cohomology is 
$\cB_H=\{q_i ,\q_i\a_1, \q_i\b_1, 
\q_i\g_1 | 1\leq i\leq 6\}$ and the basis given by the quantum
cohomology is $\cB_{QH}=\{q_i ,\q_i\a, \q_i\b, \q_i\g| 1\leq i\leq 6\}$.
We must have
$$
\left\{\begin{array}{ll}
  \q_i &=\q_i \\
  \q_i\a&=\q_i\a_1 \\
  \q_i\b&=\q_i\b_1 +A_1 \q_i \\
  \q_i\g&=\q_i\g_1 +A_2 \q_i\a_1
\end{array}\right.
$$
for $1\leq i\leq 6$ and constants $A_1$, $A_2$. We may compute these constants as
above. First we need the intersection pairings
$$
  \la\q_1\cup\q_4\a_3\ra_{\M}=-4, \qquad 
  \la\q_1\cup\q_4\a_1\b_1\ra_{\M}=4,\qquad 
  \la\q_1\cup\q_4\g_1\ra_{\M}=-4.
$$
The restriction of the cohomology classes involved to $N$ are as follows
$$
\left\{\begin{array}{l}
  \q_i=-\p_i h \\
  \q_i\a_1=-4\p_i\o h -\p_ih^2 \\
  \q_i\b_1 =-\p_i h^3 =4\p_i\o h^2+8\p_i\o^2 h\\
  \q_i\g_1= 2\p_i\o h^3= -8\p_i\o^2 h^2
\end{array}\right.
$$
And the Gromov-Witten invariants of $N$ that we need are
$$
\Q^N_l(\p_1z_1,\p_4z_2,z_3)=\left\{\begin{array}{ll}
  {(-8)^{b-5}\over (b-5)!}\p_1\p_2\o^2=2 {(-8)^{b-5}\over (b-5)!}, \qquad &
  b\geq 5 \\ 0, & b<5 \end{array} \right.
$$
for $z_i=\o^{a_i} h^{b_i} \in H^*(N)$, $1\leq i\leq 3$, $0\leq a_i
\leq 3$, $0\leq b_i\leq 2$, $a=a_1+a_2+a_3$, $b=b_1+b_2+b_3$ and $a+b=7$.
Then we get
$$
\left.
\begin{array}{l}
\Q_A^{\M}(\q_1,\b,\q_4\g_1)=\la \q_1\b,\q_4\g_1 \ra =
\la A_1\q_1,\q_4\g_1 \ra=-4 A_1 \\
\Q_l^N(-\p_1h,h^2,-8\p_4\o^2 h^2)=16
\end{array} \right\} \implies A_1=-4
$$
Analogously $A_2=-4$.

Regarding to the piece corresponding to the representation 
$\L^2_0 H_3$, we know from the proof of~\cite[lemma 14]{quantum} 
that $\q_i\q_j=\q_i\cup\q_j$, for any $1\leq i,j\leq 6$.
We have thus proven the following

\begin{thm}
\label{thm:main}
For the Riemann surface $\S$ of genus $g=3$, the isomorphism $H^*(\M)\isom 
QH^*(\M)$ is given by
$$
\left\{\begin{array}{lclclclcll}
\a^2 &=&\a_2  \\
\a\b &=&\a_1\b_1 &+&4 \a  \\
\b^2 &=&\b_2 &+&16 \b -8\a_2  \\
\a\g &=&\a_1\g_1 &-&3\b-3 \a_2 &-&1 \\
\b\g &=&\b_1\g_1&-&20\g-12\a_1\b_1 &+&24 \a \\
\g^2 &=&\g_2 &+& 8\g_1\a_1-6\b_2 &-&24 \a_2 -\b&-&8 \qquad\\
\q_i\a&=&\q_i\a_1, &&&&&& 1\leq i\leq 6\\
\q_i\b&=&\q_i\b_1 &-&4 \q_i, &&&&& 1\leq i\leq 6\\
\q_i\g&=&\q_i\g_1 &-&4 \q_i\a_1, &&&&&1\leq i\leq 6\\
\q_i\q_j &=\q_i\cup\q_j, &&&&&&& 1\leq i,j\leq 6
\end{array}\right.
$$ 
\end{thm}

\section{Some computations of Gromov-Witten invariants of degree two}
\label{sec:degree-2}

To finish it only remains to compute the following Gromov-Witten 
invariants of degree two: 
$\Q_{2A}^{\M} (\b,\b,\point)$ and $\Q_{2A}^{\M} (\a,\g,\point)$.
Fix three points $P_0,P_1,P_2\in \PP^1$. Then
$\Q_{2A}^{\M} (\a,\g,\point)$ equals
$$
  \#\{f:\PP^1 \ar \M | 
  \text{$f$ is holomorphic}, f_*[\PP^1]=2A,
  f(P_0)=\point, f(P_1) \in V_{\a}, f(P_2) \in V_{\g} \},
$$
and analogously for the other. Fix a generic point $\point\in \M$, 
corresponding to a bundle $E$ on $\S$. 
We recall the following result from~\cite[section 6]{Sam} (here
we suppose that the degree of the determinant $\L=\det E$ is one).

\begin{prop}
\label{prop:Sam}
  Let $\cU \ar \S\x\M$ be the universal bundle. For any holomorphic map
  $f:\PP^1\ar \M$ such that $f_*[\PP^1]=2A$ put $\EE=(1\x f)^*\cU\ar 
  \S\x\PP^1$.
  Denote by $p$ and $q$ the projections of $\S\x\PP^1$ onto $\S$ and 
  $\PP^1$ 
  respectively. Then one and only one of the following two cases hold
\begin{itemize}
  \item There is a stable bundle $V$ on $\S$ of rank $2$ and
  degree zero and some $x\in\S$ such that $0\ar p^* V\ox q^* \cO(1)
  \ar \EE \ar p^* \CC(x) \ar 0$.
  \item There is a line bundle $L$ on $\S$ of degree zero such
  that  $0\ar p^* L\ox q^* \cO(2) \ar \EE\ar p^* L^{-1}\ox\L \ar 0$.
\end{itemize}
\end{prop}

\begin{prop}
\label{prop:cR}
  Let 
 $$
 \cR =\{f:\PP^1 \ar \M|\text{$f$ is holomorphic}, f_*[\PP^1]=2A, 
  \point \in f(\PP^1)\}.
 $$
  endowed with the action of $\PGL$ given by reparametrization of 
  $\PP^1$. Then $\cR$ is a principal fibre bundle
  $\PGL \ar \cR \ar \PP(E)$. So $\cR$ is 
  smooth, it has dimension $\dim_\CC\cR=5$ and $\PP(E)$ parametrizes 
  ``non-parametrized'' $\PP^1 \ar \M$ representing the homology class 
  $2A$ and passing through $\point$.
\end{prop}

\begin{pf}
  The second case of proposition~\ref{prop:Sam} corresponds to the maps 
  $f: \PP^1\ar \M$ whose image lies in the subspace of bundles defined 
  by $N$ described in section~\ref{sec:4}, which is of dimension $5$.
  Therefore for a generic point, $\point\not\in N$ and hence any 
  $f\in \cR$ is in the first case of proposition~\ref{prop:Sam}.

  Now let $f\in \cR$ and put $\EE=(1\x f)^*\cU$. There exists a stable 
  bundle $V$ on $\S$ of rank $2$ and degree zero and some $x\in\S$ 
  such that 
\begin{equation}
  0\ar p^* V\ox q^* \cO(1) \ar \EE \ar p^* \CC(x) \ar 0.
\label{eqn:ext}
\end{equation}
  Let $P\in \PP^1$ with $f(P)=\point$, then restricting to $\S\x \{P\}$
  we get that $0\ar V \ar E \ar \CC(x)\ar 0$. So the possible $V$ 
  appearing for $f\in \cR$ are given by quotients of $E$ onto a 
  skyscraper sheaf supported on a single point of $\S$. This 
  corresponds to the $2$-dimensional space parametrized by $\PP(E)$.

  Moreover the exact sequence~\eqref{eqn:ext} is an elementary
  transformation which can be reversed as
\begin{equation}
  0\ar \EE \ar p^* (V\ox\cO(x))\ox q^* \cO(1) \ar p^* 
  (\CC(x)\ox\cO(x)) \ox q^*\cO(2) \ar 0.
\label{eqn:ext-dual}
\end{equation}
  So the possible bundles $\EE$ over $\S\x \PP^1$ appearing in an
  exact sequence~\eqref{eqn:ext} are parametrized by epimorphisms
  $p^* (V\ox\cO(x))\ox q^* \cO(1) \surj p^*(\CC(x)\ox\cO(x)) \ox q^*\cO(2)$, 
  i.e.\ (restricting to $\{ x \}\x \PP^1\subset \X$)
  by maps $\cO(1)\oplus \cO(1)\surj \cO(2)$, or equivalently,
  by two linear forms in $H^0(\PP^1,\cO(1))$ with no common zero. 
  Now note that $\PGL$ acts freely and transitively on these elements, 
  so the possible diagrams~\eqref{eqn:ext-dual} with $V$ fixed
  are parametrized by $\PGL$. The statement now follows easily.
\end{pf}

The image of $\cR$ in $\M$, $\{f(P)| f\in \cR, P\in \PP^1\}$, is 
described as follows. Put $R=\cR \x_{\PGL} \PP^1$, so there is a 
fibration $\PP^1 \ar R \ar \PP(E)$ and a natural map $R\ar \M$
whose image is $\{f(P)| f\in \cR, P\in \PP^1\}$. Now $R$ is smooth, 
compact and of dimension $3$. Let us denote $\a_R$, $\b_R$ and $\g_R$ 
for the pull-back of $\a$, $\b$ and $\g$ from $\M$ to $R$, 
respectively. We have the following

\begin{lem}
\label{lem:abgR}
  Let $\bar f\in H^2(\PP(E))$ denote the class of the fibre of 
  $\PP(E)\ar \S$. Then 
  $\la\g_R,[R]\ra=-12$ and $\b_R/[\PP^1]=-\bar f$,
  where $\PP^1$ is the fibre of $R\ar \PP(E)$.
\end{lem}

\begin{pf}
  The cohomology ring of $\PP(E)$ is 
  $H^*(\PP(E))=H^*(\S)[h]/(h^2-\L h)$, where
  $h$ is the hyperplane class of $\PP(E)$, since $c_1(E)=\L$.
  Now we construct the universal bundle $\cV$ parametrizing the 
  bundles $V$ which are quotients of $E$ onto a skyscraper sheaf 
  supported at one point of $\S$. We have an exact sequence
$$
  0\ar \cV \ar \pi_1^* E \ox \cO_{\PP(E)}(-1) \ar \cO_{\D}\ar 0
$$
  on $\S\x\PP(E)$, where $\pi_1: \S\x\PP(E) \ar \S$ is the projection 
  and $\D\subset \S \x\S$ is the diagonal divisor (we omit some 
  pull-backs when there is no danger of confussion). The total Chern 
  class of $\cV$ is
$$
  c(\cV)=(1-\D)(1+\pi_1^*\L-2h)=1+(\pi_1^*\L-2h-\D)+(2h-\pi_1^*\L)\D
$$

  Let us construct the space $R$. The fibre of $\cR\ar\PP(E)$ 
  over $V$ is given by the isomorphism $\Hom (V\ox\cO(x),
  \CC(x)\ox\cO(x))\isom H^0(\PP^1,\cO(1))$.
  The action of $\PGL$ is trivial on the first space and standard on 
  the second. So fibre of $R \ar \PP(E)$ over $V$ is 
  $\PP(\Hom(V,\CC(x))^{\vee})=\PP(V_x)$, where
  $V_x$ is the fibre of $V$ at $x$. Thus
\begin{equation}
  R= \PP(\pi_* (\cV\big|_{\D})),
\label{eqn:3}
\end{equation}
  for $\pi: \S \x \PP(E) \ar \PP(E)$. The Chern character of 
  $\pi_* (\cV\big|_{\D})$ is by Riemann-Roch
\begin{eqnarray*}
  \ch(\pi_* (\cV\big|_{\D}))& =& \ch\left( \pi_! (\cV\big|_{\D}) \right)=
  \pi_* \left( \ch (\cV\big|_{\D})\, \Todd T_{\S}\right) =\\
  & =& \pi_* \left( (2\D+ (\pi_1^*\L-2h-2\D)\D - \frac12 
  (\pi_1^*\L-2h)\D^2 ) (1-\frac12 K)\right) =\\
  & =& 2+ \L+K-2h 
\end{eqnarray*}
  on $\PP(E)$, where $K$ is the canonical class of $\S$ and
  using that $\pi_* \D^2=-K$. 
  We deduce that the cohomology ring of $R$ is 
  $H^*(R)=H^*(\PP(E))[k]/(k^2=c_1 k -c_2)$, where
  $k$ is the hyperplane class of~\eqref{eqn:3} and
\begin{eqnarray*}
  c_1 &=& \L +K-2h \\
  c_2&=& -2 h K
\end{eqnarray*}
  The universal bundle over $\S\x R$ is given by 
$$
  0\ar \cE \ar \cV \ox\cO(\D) \cO_{R}(-1) \ar \cO_{\D}\ox\cO(\D) 
   \ar 0,
$$
  from where we get that $c_1(\cE)=\pi^*_1\L-2h-2k$ 
  and $c_2(\cE)=(2h-\pi^*_1\L+k)\D +2hK+kK$. 
  Let $\frg_{\cE}$ be the associated $SO(3)$-bundle to $\cE$.
  We have
$$
  p_1(\frg_{\cE})=4c_2(\cE)-c_1(\cE)^2=4\D (2h+k-\pi_1^*\L).
$$
  Then $\mu_R: H_*(\S) \ar H^{4-*}(R)$ is $\mu_R(z)=-\frac14 
  p_1(\frg_{\cE}) /z$, so that
$$
\left\{\begin{array}{l} \a_R=2\mu_R(\S)=2\L -4h-2k\\
\b_R=-4\mu_R(x)= \\ -(8h+k) f
\g_R=-2\sum_{i=1}^3\mu_R(\g_i)\mu_R(\g_{3+i})=-6 (2h+k)^2 f=
  - 12 h k f 
\end{array}\right.
$$
  where $f\in H^2(R)$ is the class of the fibre of $R\ar \S$.
  Therefore $\la\g_R,[R]\ra=-12$ and $\b_R/[\PP^1]=-\bar f$ where
  $\bar f\in H^2(\PP(E))$ is the class of the fibre of
  $\PP(E)\ar \S$.
\end{pf}

Now we may finish our computation. We have
\begin{equation}
\label{eqn:ag}
 \Q_{2A}^{\M} (\a,\g,\point) =\a[2A] \la \g_R,[R] \ra=-24,
\end{equation}
since any line in $\cR$ cuts $V_{\a}$ in $\a[2A]$ points. Also
$\b_R/[\PP^1]\in H^2(\PP(E))$ is cohomology class dual to the 
subset of ``non-parametrized'' $\PP^1\ar \M$ (representing $2A$
and passing through $\point$) which lie in (a cycle 
Poincar\'e-dual to) $\b_R$. Therefore
\begin{equation}
\label{eqn:bb}
 \Q_{2A}^{\M} (\b,\b,\point)= \la (\b_R/[\PP^1])\cup 
  (\b_R/[\PP^1]) ,[\PP(E)]\ra= 0
\end{equation}

\end{document}